\documentclass[a4paper,12pt]{article}
\usepackage{amsfonts,amssymb, pdfpages,xcolor}

\title{Irregular shock waves formation as continuation of  analytic solutions}
\author{M. Colombeau,\\  \textit{mcolombeau@ime.usp.br}\\
Instituto de Matem\'atica e Estat\'istica,\\Universidade de S\~ao Paulo, Brazil}
\date{}

\begin{document}
\maketitle
\textit{Dedicated to the Memory of Vladimir Shelkovich.}\\

\begin {abstract} This paper is devoted to the blow-up of analytic solutions with the emergence of irregular solutions. We first consider the Panov-Shelkovich system which is, as far as the author knows, the first system where such types of $\delta,\delta'$-wave solutions have been explicitely exhibited.\ \  \  \

 We propose a method to reduce this system of nonlinear PDEs to a system of two ODEs in Banach spaces, which permits to obtain theoretical existence of approximate solutions for the Cauchy problem by constructing a weak asymptotic method issued from [12] and Maslov asymptotic analysis. Further this method allows to study these PDEs through their ODEs representation by basic elementary numerical schemes very easy to construct. We  then observe the expected results from the previous theoretical method; this also  gives confidence in the mathematical proofs. We  prove in scales of Banach spaces that this method gives back the  classical analytic solutions. We also study solutions in the form of $(\delta^n)'$ for other similar systems. Indeed we show formation of very irregular shock waves when the existence time of a classical analytic solution is over.\ \  \ \  \

Then we study extensions in 2-D of similar systems  exhibiting irregular solutions for which this method permits to obtain  existence of approximate solutions for the Cauchy problem. We finally sketch adaptations to provide a weak asymptotic method to the 3-D Euler-Poisson equations with application to pressureless fluid dynamics and cosmology \cite{ColombeauSiam, ColombeauNMPDE,Colombeaugravitation}.\\

\end{abstract}

\textit{This research was supported by the Funda\c c\~ao de Amparo a Pesquisa do Estado de S\~ao Paulo, processo 2012/15780-9.}\\

\textbf{1. Introduction}. Systems of conservation laws can admit irregular solutions that, further, can emerge
from regular ones: for instance $\delta$-shock waves have been intensively studied in recent years,
\cite{AlbeDani, Albeverio, AlbeShelk, Mitrovic, Shelkovich3, Mitrovic3,  ShelkovichRMS, Shelkovich1}, among other  references. $\delta$-shock wave
generation from continuous data has been studied in \cite {Mitrovic, Mitrovic2, Mitrovic3} by a method of new
characteristics that permits to understand geometrically their formation and their interaction. The $\delta^{(n)}$-shock waves were introduced and studied in
\cite{Panov, ShelkovichRMS, Shelkovichmat}. Like the $\delta-$waves these $\delta^{(n)}$-waves are put in
evidence and studied through the concept of weak asymptotic method, issued from  \cite{Danilov1}, inspired by Maslov asymptotic analysis
\cite{Maslov}, which consists in exhibiting a family of approximate
solutions $S(x,y,z,t,\epsilon)$ that tend to satisfy the equations in the sense of distributions in $x,y,z$ when the parameter $\epsilon$ tends to 0 \cite {Mitrovic, Danilov1,
Shelkovich2, Shelkovich3, Omel'yanov, Panov, ShelkovichRMS, Shelkovichmat, Shelkovich1}.\\

Motivated by the problem  to put in evidence generation
of $\delta^{(n)}$-shock waves as time continuations of analytic solutions  of the Panov-Shelkovich system, we construct a weak asymptotic method for the
Cauchy problem to a family of systems containing the original Panov-Shelkovich system, that also applies for systems in nontriangular form.\\

 In this paper we first
consider  the triangular system

\begin{equation}\frac{\partial}{\partial t}u+\frac{\partial}{\partial x}(au^2)=0,  \end{equation}
\begin{equation}\frac{\partial}{\partial t}v +\frac{\partial}{\partial x}(buv)=0,  \end{equation}  
\begin{equation} \frac{\partial}{\partial t}w+\frac{\partial}{\partial x}(cuw)+\frac{\partial}{\partial x}(P(v))=0, \end{equation}
\\
 on the torus $\mathbb{T}=\mathbb{R}/(2\pi\mathbb{Z})$,with $a,b,c\in \mathbb{R}, \ P $ a polynomial in one variable with real coefficients and arbitrary degree $n$.  The periodicity
assumption does not bring any restriction to the problem as long as one studies solutions in finite space and finite
time intervals and as long as waves propagate with finite speed, which is the case in the present paper.
With $a=1, b=2, c=2, P(v)=2v^2$ one recovers the original Panov-Shelkovich system
\cite{Panov,ShelkovichRMS,Shelkovichmat} for $\delta'$-shock waves. With $P(v)=v^n, n>2$, one obtains far more
irregular  types of solutions than $\delta^{(n)}$-shock waves. They are defined as usual through the weak asymptotic
method, see \cite{DanilovO1, DanilovO2, MaslovOmel} for similar objects such as infinitely narrow solitons.
The weak asymptotic method presented here also extends to nontriangular systems
with arbitrary $L^1_{loc}$  initial data on the multidimensional torus $\mathbb{T}^d=\mathbb{R}^d
/(2\pi\mathbb{Z})^d$. \\

We construct this weak asymptotic method by replacing the system (1-3) by a family of two ODEs  in a classical Banach
space of continuous functions. To visualize  a solution it suffices to use classical numerical approximations for the solution of these
ODEs: here we use  the explicit Euler order one  scheme in time  to observe very conveniently the creation of these irregular shock waves as time continuations of analytic Cauchy-Kovalevska
solutions and their behavior, in a rigorous mathematical way, from classical convergence results  for ODEs. \\

The method presented here can be extended to  weak asymptotic methods based on suitable a-priori estimates  for
 3-D systems of  pressureless fluids without self-gravitation  and with self-gravitation \cite{Colombeaugravitation}. The method of construction of sequences for the asymptotic analysis has been inspired from the sequences constructed for 
the  numerical schemes in \cite{ColombeauSiam, ColombeauNMPDE, Colombeauideal} which has been a preliminary step.\\

\textbf{2. Construction of a weak asymptotic method for system (1-3)}.
 We assume the
initial conditions $u_0,v_0,w_0$ are locally integrable  on the torus $\mathbb{T}=\mathbb{R}
/(2\pi\mathbb{Z})$. For simplification  we state the proof in the case $a=b=c=1$ and $P(v)=v^n, \ n \in \mathbb{N}$ since it simplifies the formulas and makes no change in the proof without any loss of generality.
The
approximate initial conditions $u_0^\epsilon,v_0^\epsilon,w_0^\epsilon, (\epsilon>0,\epsilon \rightarrow 0$),
which are regularizations of $u_0,v_0,w_0$, will be chosen periodic within this period. Therefore the approximate
solutions $u^\epsilon,v^\epsilon,w^\epsilon$ that we will construct will be periodic as well.\\

\textbf{Construction of the family ($u^\epsilon$)}.  We start from a family $(u^\epsilon)$ of approximate solutions of (1) such that
\begin{equation} \exists M_1>0 \ / \ |u^\epsilon(x,t)|\leq M_1 \  \forall \epsilon, \forall x\in \mathbb{T}, \forall t\in
[0,+\infty[. \end{equation}
For the  proof  of lemma 2 we need a family $(u^\epsilon)$ such that 
\begin{equation} \forall \beta>0 , \ \forall \delta>0 \ \exists \ const \ / \ |\frac{\partial}{\partial
x}u^\epsilon(x,t)| \leq \frac{const}{\epsilon^\beta} \ \forall x \in \mathbb{T}, \ \forall t\in
[0,\delta].\end{equation}
This property can be obtained by starting from a family $(u^\nu)$ of $2\pi$-periodic viscous solutions of (1) with viscosity
coefficient $\nu$ and by reindexing this family by means of a function $\epsilon\longmapsto
\nu(\epsilon)$ so as to produce a growth of $sup_{x,t\in [0,\delta]} |\frac{\partial}{\partial
x}u^\epsilon(x,t)|$ as slow as wanted when $\epsilon\rightarrow 0$, and for all given $\delta$. \\


\textbf{Construction of the family ($v^\epsilon$)}. To obtain the approximate solutions $v^\epsilon$ of
equation (2) we do as follows. We choose as initial condition a family $(v_0^\epsilon)$ such that
$\|v_0^\epsilon\|_{L^1(-\pi,+\pi)} \leq \ const$ independent on $\epsilon$ since $\|v_0\|_{L^1(-\pi,+\pi)} <
+\infty$. In a first step, for each $\epsilon$ we obtain existence and uniqueness of a solution of  a linear ODE in the Banach space
$(\mathcal{C}_b(\mathbb{R}), \| .\|_\infty)$ of all bounded continuous ($2\pi-$periodic) functions on $\mathbb{R}$ equipped with
the sup. norm. In a second step we obtain $v^\epsilon$ as a suitable regularization by convolution of the
solution of this ODE.\\

 To this end we note as usual 

\begin{equation}u^+(x)=max(0,u(x)), \ u^-(x)=max(0,-u(x)),\end{equation}
then
\begin{equation}u(x)=u^+(x)- u^-(x), \ |u(x)|=u^+(x)+u^-(x).\end{equation}
For $0<\epsilon<1 $ we consider the homogeneous linear ODE
\begin{equation}
\frac{d}{dt}X^\epsilon(x,t)=\frac{1}{\epsilon}[X^\epsilon(x-\epsilon,t)u^{\epsilon+}(x-\epsilon,t)-X^\epsilon(x,t)|u^{\epsilon}(x,t)|+X^\epsilon(x+\epsilon,t)u^{\epsilon-}(x+\epsilon,t)],\end{equation}
\begin{equation} X^\epsilon(x,0)=v_0^\epsilon(x).\end{equation}
\\
\textit{Remark}. Formula (8) is an approximation of (2), i.e. $\frac{\partial}{\partial t}X+\frac{\partial}{\partial x}(Xu)=0$. 
In the particular case $u^\epsilon(x,t)$, refered to as a "velocity",  has a constant sign, the
right-hand side of formula (8) reduces at once to the classical discretization of the $x-$derivative $\frac{\partial}{\partial x}(Xu)$
to the left for positive velocity $u$ and to the right for negative velocity. This formula that replaces the
$x$-derivative is issued from the scheme in \cite{ColombeauSiam,ColombeauNMPDE,Colombeauideal} and therefore it
can be obtained intuitively as follows considering $X$ as a mass density and $u$ as a velocity of matter. Consider cells of length $\epsilon$ centered respectively at the points
$x,x-\epsilon,x+\epsilon$. Then we state that the quantity $\epsilon.X^\epsilon(x,t+dt)$ in the cell of center
$x$ at time $t+dt$ is equal to the quantity $\epsilon.X^\epsilon(x,t)$ at time $t$ minus the quantity
$(X^\epsilon |u^\epsilon|.dt)(x,t)$ escaped from this cell between $t$ and $t+dt$ (to the left if
$u^\epsilon<0$, to the right if $u^\epsilon>0$) plus the quantities received by this cell from the left and from the
right, namely $(X^\epsilon .u^{\epsilon +}.dt) (x-\epsilon,t), (X^\epsilon .u^{\epsilon -}.dt) (x+\epsilon,t)$
respectively, under the assumption $|u^\epsilon dt| \leq \epsilon$ so that all matter comes from the two
closest neighbor cells only. One lets $dt$ tend to 0 for fixed $\epsilon$ to obtain the differential equation
(8) instead of the numerical scheme in \cite{ColombeauSiam,ColombeauNMPDE,Colombeauideal} in which the "infinitesimal $dt$" is replaced by a fixed value $\Delta t$ imposing from the requirement $ \|u\|_\infty\Delta t\leq \epsilon$ a bound on the velocity. This formulation of
the derivative $\frac{\partial}{\partial x}(Xu)$ is more precise than a classical quotient in that it takes into account the real
phenomenon of transport of the density $X$ by the velocity $u$ whatever the (variable) sign of $u$ is. Further the
formulation (8) permits to use as mathematical tool for fixed $\epsilon>0$ the theory of ordinary differential
equations in Banach spaces.\\


Since, from (4), the coefficients $u_\epsilon^\pm$ are in $ \mathcal{C}_b(\mathbb{R})$, the linear ODE (8, 9)
has a global unique solution $X^\epsilon$ which is a continuous function on $[0,+\infty[$ valued in the Banach
space $(\mathcal{C}_b(\mathbb{R}), \| .\|_\infty)$. Then one sets
\begin{equation} v^\epsilon(x,t)=(X^\epsilon(.,t)*\phi_{\epsilon^\alpha})(x)\end{equation}
where $\phi$ is a $\mathcal{C}^\infty $ function on $\mathbb{R} $ with compact support such that
$\int\phi(x)dx=1$, $\phi_{\epsilon^\alpha}(x)=\frac{1}{\epsilon^\alpha}\phi(\frac{x}{\epsilon^\alpha}),
\alpha\in ]0,1] $ to be chosen later. Note that $v^\epsilon$ is of class $\mathcal{C}^\infty$ in $x$ while
$X^\epsilon$ is only $\mathcal{C}^0$ in $x$.\\

\textbf{Lemma 1.}
\begin{equation}\exists C>0 \ / \ \ \forall \epsilon>0 \ \int_{-\pi}^{+\pi}|X^\epsilon(x,t)| dx\leq \int_{-\pi}^{+\pi}|v_0^\epsilon(x)|dx\leq C. \end{equation}
\\
\textit{proof}. In case $X$ would be positive this would follow at once from an integration of (8) on a period using (7).
But $X$ can have arbitrary sign therefore a different proof is needed. Developping $X^\epsilon(x,t+dt)$ from (8) according to Taylor's formula we obtain\\
\\
$X^\epsilon(x,t+dt)=\frac{dt}{\epsilon} X^\epsilon (x-\epsilon,t)u^{\epsilon
+}(x-\epsilon,t)+(1-\frac{dt}{\epsilon}|u^\epsilon(x,t)|)X^\epsilon (x,t) +
\frac{dt}{\epsilon} X^\epsilon (x+\epsilon,t) u^{\epsilon -}(x+\epsilon,t)+dt.o^\epsilon (x,t)(dt)$\\  
\\
where $o^\epsilon(x,t)(dt)\rightarrow 0$ uniformly in $(x,t), \ x\in \mathbb{R}, \ t\in I\subset\subset [0,+\infty[$, when $dt\rightarrow 0$. 
This uniform bound on $o^\epsilon(x,t)(dt)$ follows from the application to the $\mathcal{C}^1 $ map
$t\longmapsto X^\epsilon(.,t)$ of the mean value theorem under the form: if $f$ is a
$\mathcal{C}^1$ function then $\|f(t+dt)-f(t)-f'(t).dt\| \leq sup_{0<\theta<1}\|(f'(t+\theta dt)-f'(t)).dt\|$ and
from the uniform continuity of $\frac{dX^\epsilon(.,t)}{dt}$ in $t$, taking values in the Banach space
$\mathcal{C}_b(\mathbb{R})$, on any compact set in $t$-variable. 
 In order to have
$1-\frac{dt}{\epsilon}|u^\epsilon(x,t)|\geq 0$, we choose  $dt>0$
small enough such that $\frac{dt M_1}{\epsilon}\leq 1$.
Therefore taking absolute values,\\

$|X^\epsilon(x,t+dt)|\leq\frac{dt}{\epsilon}|X^\epsilon(x-\epsilon,t)|u^{\epsilon
+}(x-\epsilon,t)+(1-\frac{dt}{\epsilon}|u^\epsilon(x,t)|)|X^\epsilon(x,t)| +
\frac{dt}{\epsilon}|X^\epsilon(x+\epsilon,t)|u^{\epsilon -}(x+\epsilon,t)+dt. |o^\epsilon (x,t)(dt)|.$\\
\\
After  changes in $x$-variable and use of periodicity for simplification of integrals\\

$\int_{-\pi}^{+\pi}|X^\epsilon(x,t+dt)|dx\leq\frac{dt}{\epsilon}\int_{-\pi}^{+\pi}|X^\epsilon(x,t)|u^{\epsilon
+}(x,t)dx+\int_{-\pi}^{+\pi}[1-\frac{dt}{\epsilon}u^{\epsilon +}(x,t)-\frac{dt}{\epsilon}u^{\epsilon
-}(x,t)]|X^\epsilon(x,t)| dx+
\frac{dt}{\epsilon}\int_{-\pi}^{+\pi}|X^\epsilon(x,t)|u^{\epsilon -}(x,t)dx+dt \int_{-\pi}^{+\pi}|o^\epsilon
(x,t)(dt)|dx=
\int_{-\pi}^{+\pi}|X^\epsilon(x,t)|dx+dt.o_1^\epsilon (t)(dt) $\\
where $o_1^\epsilon (t)(dt)\rightarrow 0$ when $dt\rightarrow 0$, uniformly when $t$ ranges in a bounded set of $[0,+\infty[$.\\

We divide the interval $[t,t+\tau]$ into subintervals $[t+i\frac{\tau}{n}, t+(i+1)\frac{\tau}{n}], \ 0\leq
i\leq n-1, n$ as large as needed, and apply the above bound with $dt=\frac{\tau}{n}$ in each subinterval, i.e.\\
\\
$$\int_{-\pi}^{+\pi}|X^\epsilon(x,t+(i+1)\frac{\tau}{n})|dx\leq
\int_{-\pi}^{+\pi}|X^\epsilon(x,t+i\frac{\tau}{n})|dx+\frac{\tau}{n}o_2^\epsilon(\frac{\tau}{n})$$\\
\\
with $o_2^\epsilon$ independent on $t$ from the uniformness $o_1^\epsilon$ in $t$ : $o_2^\epsilon(dt)=sup_t o_1^\epsilon (t)(dt) $. Adding on all
$i$ one obtains \\
\\
$$\int_{-\pi}^{+\pi}|X^\epsilon(x,t+\tau)|dx\leq \int_{-\pi}^{+\pi}|X^\epsilon(x,t)|dx+n
\frac{\tau}{n}o_2^\epsilon(\frac{\tau}{n}).$$\\
\\
Since $o_2^\epsilon(\frac{\tau}{n})\rightarrow 0$ when $n\rightarrow +\infty$ we have \\
$$\int_{-\pi}^{+\pi}|X^\epsilon(x,t+\tau)|dx\leq \int_{-\pi}^{+\pi}|X^\epsilon(x,t)|dx.$$
\\
Finally, replacing $t$ by $0$ and $\tau$ by $t$ one obtains
$$\int_{-\pi}^{+\pi}|X^\epsilon(x,t)|dx\leq
\int_{-\pi}^{+\pi}|X^\epsilon(x,0)|dx=\int_{-\pi}^{+\pi}|v_0^\epsilon(x)|dx\leq const.$$ $\Box$
 
\textbf{Corollary}. \textit{There is a constant $M_2>0$, independent on $\epsilon$, such that $\forall x\in
\mathbb{R}, \forall t\geq 0$}

\begin{equation} |v^\epsilon(x,t)|\leq \frac{M_2}{\epsilon^\alpha}, \ |\frac{\partial}{\partial x}v^\epsilon(x,t)|\leq
\frac{M_2}{\epsilon^{2\alpha}}. \end{equation}
\\
\textit{proof}. From (10) $|v^\epsilon(x,t)|=|\int
X^\epsilon(x-y,t)\frac{1}{\epsilon^\alpha}\phi(\frac{y}{\epsilon^\alpha})dy|\leq
\int_{supp \ \phi}|X^\epsilon(x-y,t)|\frac{1}{\epsilon^\alpha}\|\phi\|_\infty dy\leq const\frac{1}{\epsilon^\alpha}
$ from lemma 1. We obtain the second formula by a similar proof.$\Box$\\

\textbf{Lemma 2}. \textit{The family $(u^\epsilon,v^\epsilon)$ is a weak asymptotic method for equation (2)
when $\epsilon \rightarrow 0$ if $0<\beta<\alpha$.}\\
\\
\textit{proof.} We have to prove that $\forall \psi  \ \mathcal{C}^\infty$ on $\mathbb{R}$ with compact support
\begin{equation} I:=\int[\frac{\partial}{\partial t}(v^\epsilon)(x,t)\psi(x) -u^\epsilon(x,t)v^\epsilon(x,t)\psi'(x)]dx\rightarrow
0\end{equation}
\\
 when $\epsilon\rightarrow 0$. From (10)
 
$$I=\int[(\frac{d}{dt}X^\epsilon(.,t)*\phi_{\epsilon^\alpha})(x)\psi(x)
-u^\epsilon(x,t)(X^\epsilon(.,t)*\phi_{\epsilon^\alpha})(x)\psi'(x)]dx.$$
 Now, from (8), since $|u|=u^++u^-$\\
 
$I=\int\{\frac{1}{\epsilon}[(X^\epsilon u^{\epsilon +})(x-\lambda-\epsilon,t)-(X^\epsilon u^{\epsilon
+})(x-\lambda,t)-(X^\epsilon u^{\epsilon -})(x-\lambda,t)+(X^\epsilon u^{\epsilon
-})(x-\lambda+\epsilon,t)]\frac{1}{\epsilon^\alpha}\phi(\frac{\lambda}{\epsilon^\alpha})\psi(x)-u^\epsilon(x,t)X^\epsilon(x-\lambda,t)\frac{1}{\epsilon^\alpha}\phi(\frac{\lambda}{\epsilon^\alpha})\psi'(x)\}dxd\lambda.$\\
 \\
 Changes in $x$-variable give\\
 
$I=\int\{\frac{1}{\epsilon}[(X^\epsilon u^{\epsilon +})(x-\lambda,t)(\psi(x+\epsilon)-\psi(x))-(X^\epsilon
u^{\epsilon
-})(x-\lambda,t)(\psi(x)-\psi(x-\epsilon))]-u^\epsilon(x,t)X^\epsilon(x-\lambda,t)\psi'(x)\}\frac{1}{\epsilon^\alpha}\phi(\frac{\lambda}{\epsilon^\alpha})dxd\lambda.$\\
\\
Let $K$ be a finite interval in $\mathbb{R}$ containing support of $\psi$ and its translation by $\pm\epsilon$.
We use that\\ 
\\
$\int (X^\epsilon u^{\epsilon +})(x-\lambda,t)\frac{\psi(x+\epsilon)-\psi(x)}{\epsilon}dx=\int_K
(X^\epsilon u^{\epsilon +})(x-\lambda,t)(\psi'(x)+O_{[x]}(\epsilon))dx$\\
where the notation $O_{[x]}(\epsilon)$ means a dependence on $x$ that disappears in the next bound. From   lemma 1 and from (4)\\
\\ 
$\int_K (|X^\epsilon| u^{\epsilon +})(x-\lambda,t)|O_{[x]}(\epsilon)|)dx\leq  M_1 |O(\epsilon)| const.$
\\
\\ 
Therefore, with another $O(\epsilon)$,\\

$I=\int[(X^\epsilon u^{\epsilon +})(x-\lambda,t)\psi'(x)-(X^\epsilon u^{\epsilon
-})(x-\lambda,t)\psi'(x)-u^\epsilon(x,t)X^\epsilon(x-\lambda,t)\psi'(x)]\frac{1}{\epsilon^\alpha}\phi(\frac{\lambda}{\epsilon^\alpha})dxd\lambda+
O(\epsilon).$\\
\\
Since $u^{\epsilon +}-u^{\epsilon -}=u^\epsilon$, and after a change of variable\\

$I=\int X^\epsilon(x-\epsilon^\alpha \mu,t)[u^\epsilon(x-\epsilon^\alpha
\mu,t)-u^\epsilon(x,t)]\psi'(x)\phi(\mu)dxd\mu+O(\epsilon).$\\
\\
From (5), application of the mean value theorem in $u^\epsilon$ gives $|u^\epsilon(x-\epsilon^\alpha
\mu,t)-u^\epsilon(x,t)|\leq \frac{const}{\epsilon^\beta}\epsilon^\alpha |\mu|$. Then from    lemma 1,\\
\\ 
$|I|\leq \|X^\epsilon\|_{L^1(K)} \frac{const}{\epsilon^\beta}\epsilon^\alpha\int|\mu \phi(\mu)|d\mu 
+O(\epsilon)\leq const.\epsilon^{\alpha-\beta}+O(\epsilon)$. $\Box$ \\
\\

\textbf{Construction of the family ($w^\epsilon$)}. Now we proceed to the construction of
$(w^\epsilon)_\epsilon$. We define $w^\epsilon$ as the solution of the following linear ODE with second
member in the Banach space $\mathcal{C}_b(\mathbb{R})$

\begin{equation} \frac{d}{dt}w^\epsilon(x,t)=\frac{1}{\epsilon}[(w^\epsilon u^{\epsilon +})(x-\epsilon,t)-
(w^\epsilon |u^{\epsilon}|)(x,t)+(w^\epsilon u^{\epsilon -})(x+\epsilon,t)]- n(v^\epsilon)^{(n-1)}(x,t)
\frac{\partial}{\partial x}v^\epsilon(x,t),\end{equation}
\\
(notice that from (10) $v^\epsilon$ is a $\mathcal{C}^\infty$ function so that the $x$-derivative above makes
sense), with initial condition \\
\begin{equation} w^\epsilon(x,0)=w_0^\epsilon(x), \end{equation}
\\
where  $w_0^\epsilon \in \mathcal{C}^\infty(\mathbb{T}), \ \    \|w_0^\epsilon-w_0\|_{L^1(\mathbb{T})} \rightarrow 0$. Since $u^\epsilon$ and $v^\epsilon$ are known, the equation (14) with unknown $w^\epsilon$ is linear with second member and it has bounded
coefficients $u^\epsilon$. Therefore, for each $\epsilon>0$ it admits a unique global $\mathcal{C}^1$ solution
$t\longmapsto w^\epsilon(t)$ valued in the Banach space $\mathcal{C}_b(\mathbb{R})$ .\\
\\
\textbf{Lemma 3}.\begin{equation} \exists C>0 \ / \ \ \forall \delta>0, \ \forall t\in[0,\delta], \ \forall \epsilon>0, \ \int_{-\pi}^{\pi}|w^\epsilon(x,t)|dx \leq
\frac{C}{\epsilon^{(n+1)\alpha}}.\end{equation}\\
\\
\textit{proof.} If $dt>0$ is small enough it follows from (14) that\\

 $w^\epsilon(x,t+dt)=
w^\epsilon(x,t)+
\frac{dt}{\epsilon}[(w^\epsilon u^{\epsilon +})(x-\epsilon,t)- (w^\epsilon |u^{\epsilon}|)(x,t)+$
\begin{equation}(w^\epsilon u^{\epsilon -})(x+\epsilon,t)]- n dt (v^\epsilon)^{(n-1)}(x,t) \frac{\partial}{\partial x}
v^\epsilon(x,t)+ dt. o^\epsilon(x,t)(dt), \end{equation}
where $o^\epsilon(x,t)(dt)$ tends to 0 when $dt\rightarrow 0$ uniformly if $x\in \mathbb{R}$ and $t$ in a
compact set in $[0,+\infty[$  (the proof is the same as the one in lemma 1).\\
Therefore, since for $dt>0$ small enough depending on $\epsilon$, \ one has from (4)
$1-\frac{dt}{\epsilon}|u^\epsilon(x,t)|\geq 0 \ \forall x\in \mathbb{R}, \forall t \geq 0$. It follows that\\
\\
$\int_{-\pi}^{+\pi}|w^\epsilon(x,t+dt)|dx \leq \int_{-\pi}^{+\pi}\frac{dt}{\epsilon}(|w^\epsilon|u^{\epsilon
+})(x-\epsilon,t)dx+\int_{-\pi}^{+\pi}(1-\frac{dt}{\epsilon}|u^\epsilon|(x,t))|w^\epsilon(x,t)|dx+\int_{-\pi}^{+\pi}\frac{dt}{\epsilon}(|w^\epsilon|u^{\epsilon
-})(x+\epsilon,t)dx+ndt\int_{-\pi}^{+\pi}|(v^\epsilon)^{n-1}| . |\frac{\partial}{\partial x}v^\epsilon|(x,t)dx +dt .
o^\epsilon(dt)$\\
\\
for $dt>0$ small enough depending on $\epsilon$. Here $o^\epsilon(dt)\leq2\pi sup_{x,t}|o^\epsilon(x,t)(dt)|$
tends to 0 uniformly ($t\in I \subset\subset [0,+\infty)$) when $dt\rightarrow 0$. Changes
in variable, $|u^\epsilon|=u^{\epsilon +}+u^{\epsilon -}$, and the bounds (12) give

$$\int_{-\pi}^{+\pi}|w^\epsilon(x,t+dt)|dx \leq $$
$\frac{dt}{\epsilon}\int_{-\pi+\epsilon}^{+\pi+\epsilon} (|w^\epsilon|u^{\epsilon +})(x,t)dx+\int_{-\pi}^{+\pi}
|w^\epsilon(x,t)|dx
     -\frac{dt}{\epsilon}\int_{-\pi}^{+\pi} |w^\epsilon(x,t)| u^{\epsilon +}(x,t) dx\\ 
-\frac{dt}{\epsilon}\int_{-\pi}^{+\pi} |w^\epsilon(x,t)|u^{\epsilon -}(x,t)dx +\frac{dt}{\epsilon}
\int_{-\pi-\epsilon}^{+\pi-\epsilon} (|w^\epsilon|u^{\epsilon -})(x,t)dx+dt\frac{const}{\epsilon^{(n+1)\alpha}}
+dt .o^\epsilon(dt).$ \\
\\
The periodicity of initial conditions, coefficients and second member implies periodicity of the solutions and
therefore simplifications. One obtains
$$\int_{-\pi}^{+\pi}|w^\epsilon(x,t+dt)|dx \leq
\int_{-\pi}^{+\pi}|w^\epsilon(x,t)|dx+dt\frac{const}{\epsilon^{(n+1)\alpha}}+dt . o^\epsilon(dt).$$
Sharing the interval $[t,t+\tau], \tau>0,$ into small subintervals $[t+i\frac{\tau}{m},t+(i+1)\frac{\tau}{m}],
\ 0\leq i\leq m-1$, applying the above inequality with $dt=\frac{\tau}{m}$, adding on $i$ and letting $m\rightarrow
\infty$ (as in lemma 1) one obtains that

$$\int_{-\pi}^{+\pi}|w^\epsilon(x,t+\tau)|dx \leq
\int_{-\pi}^{+\pi}|w^\epsilon(x,t)|dx+\tau\frac{const}{\epsilon^{(n+1)\alpha}}.$$
Finally setting $t=0$ and $\tau=t$  one obtains (16).$\Box$\\

In conclusion of this construction let us recall the assumptions done throughout it. \\

$\bullet$  First, all Cauchy data $u_0,v_0,w_0$ are periodic and their regularizations
$u_0^\epsilon,v_0^\epsilon,w_0^\epsilon$ are chosen also periodic with same period, i.e. we consider the
problem on the one dimensional torus.\\

$\bullet$   As usual $u_0$ is $L^\infty$ and one considers approximate solutions $u^\epsilon$ of equation (1) satisfying
(5): a  reindexation in $\epsilon$ of the classical viscous solutions.\\

$\bullet$   The initial conditions $v_0,w_0$ are $L^1_{loc}$, and their regularizations $v_0^\epsilon,w_0^\epsilon$ are
also chosen  $L^1_{loc}$ uniformly in $\epsilon$.\\

Then we are going to prove:\\
\\
\textbf{Theorem 1. Provided $0<\alpha<\frac{1}{n+1}$ and $0<\beta<\alpha$ the family} $(u^\epsilon, v^\epsilon,
w^\epsilon)_\epsilon$ \textbf{ is a weak asymptotic method for system (1, 2, 3).}\\
\\
\textit{Proof.} From the choice of classical approximate solutions for equation (1) and from lemma 2 for
equation (2) it remains to prove that $\forall \psi \in \mathcal{C}^\infty(\mathbb{R})$ with compact support\\
$$J:=\int[\frac{\partial}{\partial t}(w^\epsilon)\psi-u^\epsilon w^\epsilon\psi'-(v^\epsilon)^n\psi']dxdt\rightarrow 0$$
when $\epsilon\rightarrow 0$. From (14) and $|u|=u^++u^-$,\\

$J=\int\{[\frac{1}{\epsilon}(w^\epsilon u^{\epsilon +})(x-\epsilon,t)-\frac{1}{\epsilon}(w^\epsilon
(u^{\epsilon +}
+u^{\epsilon -}))(x,t)+\frac{1}{\epsilon}(w^\epsilon u^{\epsilon -})(x+\epsilon,t) -\frac{\partial}{\partial x} ((v^\epsilon)^n)
(x,t)]\psi(x)-(u^\epsilon w^\epsilon)(x,t)\psi'(x) -(v^\epsilon)^n(x,t)\psi'(x)\}dxdt.$\\
\\
The two terms involving $(v^\epsilon)^n$ simplify each other. After changes in $x$-variable\\

$J=\int\{[\frac{1}{\epsilon}(w^\epsilon u^{\epsilon +})(x,t)\psi(x+\epsilon)-\frac{1}{\epsilon}(w^\epsilon
u^{\epsilon +})(x,t)\psi(x) -\frac{1}{\epsilon}(w^\epsilon u^{\epsilon
-})(x,t)\psi(x)+\frac{1}{\epsilon}(w^\epsilon u^{\epsilon -})(x,t)\psi(x-\epsilon)
-(u^\epsilon w^\epsilon)(x,t)\psi'(x)\}dxdt.$\\
\\
If $I$ is a finite interval containing the support of $\psi$ and the translated of this support by
$\pm\epsilon$ one has\\

$\int_I ( w^\epsilon u^{\epsilon \pm})(x,t) \frac{\psi(x+\epsilon)-\psi(x)}{\epsilon}dx =\int_I (w^\epsilon
u^{\epsilon \pm}) (x,t) \psi'(x)dx+\int_Iw^\epsilon u^{\epsilon \pm} (x,t)O_{[x]}(\epsilon) dx$ \\
\\
and,  from (4, 16) \\
\\
$|\int_I (w^\epsilon u^{\epsilon \pm}) (x,t)O_{[x]}(\epsilon) dx| \leq M_1 \frac{const}{\epsilon^{(n+1)\alpha}}
\epsilon=const .\epsilon^{(1-(n+1)\alpha)}.$\\

We come back to $J$. Using the above and $u^\epsilon=u^{\epsilon +}-u^{\epsilon -}$, the terms involving $\psi'$ disappear. Choosing  $\alpha<\frac{1}{n+1}$ and $\beta<\alpha$ to
apply lemma 2, the quantity $J$ tends to 0 when $\epsilon \rightarrow 0$. $\Box$\\

\textbf{3. Numerical confirmations}. The theoretical weak asymptotic method introduced in this paper permits to reduce the study of approximate solutions to the Cauchy problem for systems of PDEs such as (1-3) to a system of ODEs. Since this process is constructive one can approximate numerically the solutions. We observed the results expected from the above proofs.
 This permits to visualize the $\delta'$-waves and brings a confirmation of the theoretical results.\\
\vskip 8 cm
\textit{figure1.  Emergence of a $\delta'$ shock wave from an analytic solution.}\\

 In the sequel we will visualize $\delta"$-waves from a corresponding Panov-Shelkovich system, and more general objects such as derivatives of powers of the Dirac measure. 
This can be done very easily from basic elementary numerical methods for ODEs. First, let us notice that the linearity of the ODEs (8, 14) implies the convergence of the explicit Euler order one method for fixed $\epsilon$.
 Therefore there is no lack of rigor in using the numerical method to approximate the ODEs.

\begin{figure}[h]
\centering
\includepdf[width=\textwidth]{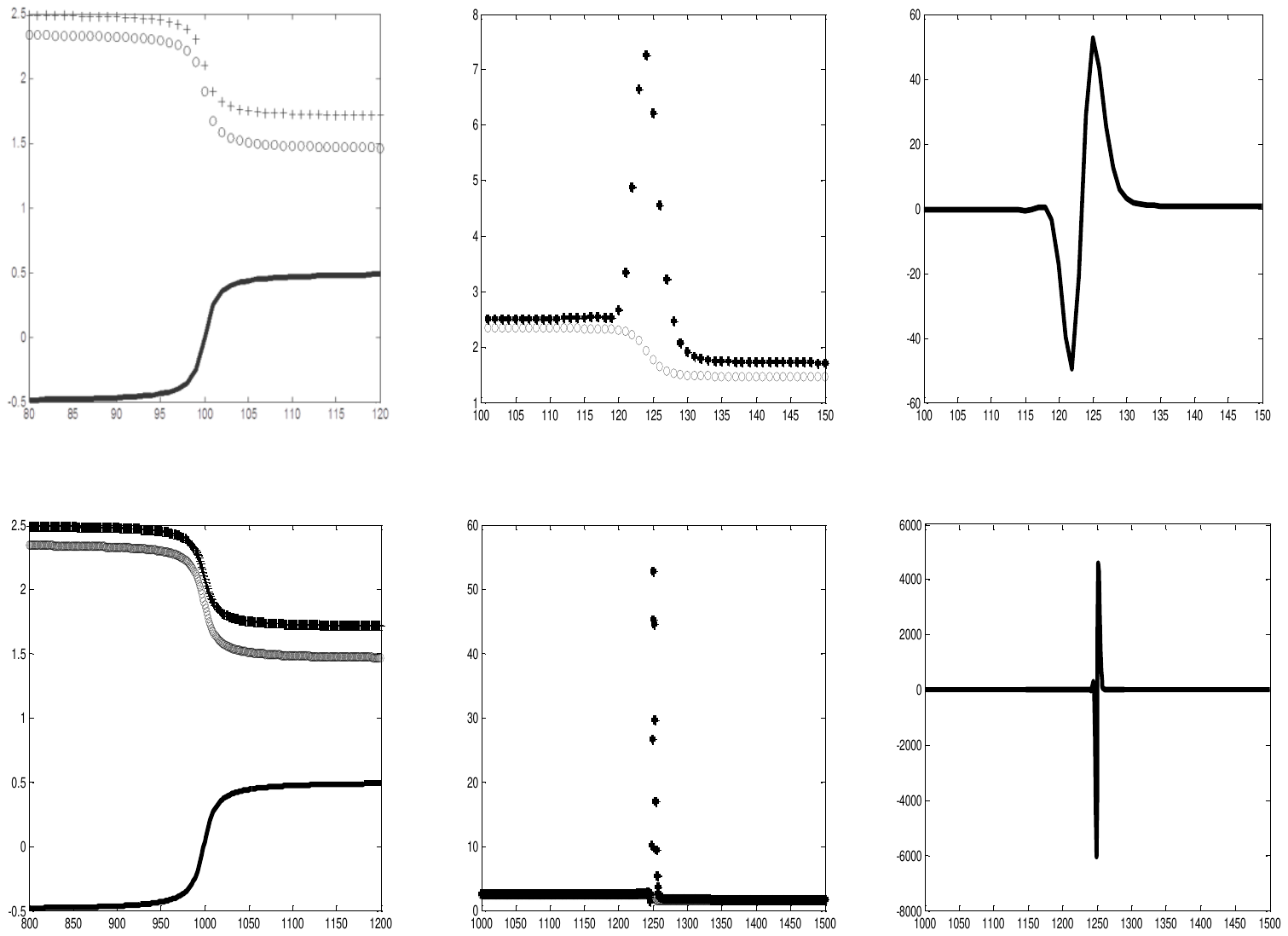}
\end{figure}

 In the tests
below it has sufficed to use the explicit order one Euler method, and to reproduce convolution  (10)  by an averaging. In numerical calculations, property (5) follows from a similar averaging as  the one in \cite{ColombeauNMPDE,Colombeauideal}. The scheme so obtained is very close to the scheme used
in \cite{ColombeauSiam,ColombeauNMPDE,Colombeauideal} for systems of fluid dynamics. The novelty here is that  one observes the
emergence of $\delta'$-shock waves in $w$ as a time continuation of the Cauchy-Kovalevska analytic solution when it
ceases to exist.\\


In figure 1 one first constructs  an analytic initial condition in $u,v,w$ at time $t=0$ (left panels). This initial condition could be a restriction to the interval under concern of a set
of analytic periodic functions with larger period because of the boundedness of the "`velocity"' $u$: in other
words the periodic assumption causes no significative loss in generality for a study of the equations in finite
space and time intervals since the velocity $u$ is bounded. One observes numerically that this analytic
solution ceases to exist at time $t = 0.037$, with the well known creation of a classical shock
wave in $u$, as well as the emergence of a $\delta$-wave in $v$ and a $\delta'$-wave in $w$. These last two
waves grow with time and the results are given at time $t=1$. The functions $u,v,w$ are respectively noted with
o,+ (bold + in middle panels) and continuous line. The top figures have been done with a large value of the
space step so as to identify more easily the solution: it gives "`thick"' shock waves. The middle and right  bottom figures have
been done with a space step ten times smaller for a better visualization of the waves. In the middle
 panels we observe the usual
shock wave in $u$ and the $\delta$-shock wave in $v$. In the right panels we observe the $\delta'$-shock wave
in $w$. Of course the aspect of the analytic initial conditions does not depend on the thickness of the space
step (left panels). \\

Other calculations have been done by replacing $v^2$ in equation (3) by $v^n$ for $ n=3,4,5,\dots$. One observes
mathematically new highly singular shock waves in $w$, not pointed out by previous authors, all stemming as time continuations of analytic solutions after their analytic blow-up. These
singular waves in $w$ are located on a point; they take considerably larger values than the $\delta'$-waves for
same $\epsilon$. This suggests that the now well known $\delta,\delta',\delta^{(n)}$ waves are only particular cases of a
general phenomenon that continues the Cauchy-Kovalevska solutions after their analytic blow-up. All these one
dimensional calculations are quasi instantaneous on any standard PC.\\

We numerically identify  some of these singular shock waves for $P(v)=v^n, n=3,4,5$. For convenience we consider the Riemann
problem $u_l=2,u_r=1,v_l=2,v_r=1, w_l=0,w_r=0$ at $x=0$ (similar results can be obtained from analytic initial
conditions as in figure 1). The segment $[-0.5, 0.5]$ is divided into cells of length $\epsilon$ and we observe
the solution $w$ at time $t=1$. The observed solution is null except on a very small region in which it has the shape
of a $\delta'$ as in figure 1, but we will observe it is not a $\delta'$ if $n\not=2$. To this end we compute a primitive of this solution : it
looks like a Dirac delta measure and we compute the area of the region between the graph of the primitive and the $x$-axis, that
should be constant (for fixed $t$: here $t=1$) as a function of $\epsilon$ in case this primitive would be a Dirac delta measure. We obtain :\\

\begin{tabular}{|c|c|c|c|c|c|c|}
  
 \hline
$ n$ & $\epsilon=10^{-3}$ & $\epsilon=5. 10^{-4}$ & $\epsilon=25. 10^{-5}$ &$\epsilon=125. 10^{-6}$
&$\epsilon=625. 10^{-7}$&$\epsilon=3125. 10^{-8}$\\
 
 \hline
 2&0.084&0.087&0.089&0.090&0.090&0.090\\
 \hline
 3&0.11& 0.23& 0.48&0.98&1.96& 4.00\\
 \hline
 4&1.36&5.56&23.2&95.7&389&1500\\
 \hline
 5&17.2&139&1168&9680&$79. 10^3$&$640. 10^3$\\
\hline 
 \end{tabular} \\
 \\     


In the case $n=2$ (i.e. Panov-Shelkovich system) the values of the areas are independent on $\epsilon$ which shows that the primitive is a Dirac delta measure, therefore $w(.,t)$
is a $\delta'$-wave. In the other cases the values of the areas  depend on $\epsilon$: for $n=3$ they are approximately
multiplied by 2 at each step to the right, and by 4, 8 if $n=4,5$ respectively. A power $\delta^\alpha$ of a
Dirac delta measure is  represented by $(\frac{1}{\epsilon}\phi(\frac{x}{\epsilon}))^\alpha$
where $\phi$ is a positive continuous functions with compact support and $\int \phi(x)dx =1$. The area between
the curve and the $x$-axis is therefore equal to  $ \int\frac{1}{\epsilon^\alpha}(\phi(\frac{x}{\epsilon}))^\alpha dx=\epsilon^{1-\alpha}\int (\phi(x))^{\alpha} dx$. For the value
$\frac{\epsilon}{2}$ the area becomes $(\frac{\epsilon}{2})^{1-\alpha}\int (\phi(x))^{\alpha} dx$. For $n=3$ we numerically observe  that 
this value is twice the value of the area relative to the space step $\epsilon$: this gives $(\frac{\epsilon}{2})^{1-\alpha}=2.\epsilon^{1-\alpha}$, i.e. $\alpha=2$; for n=4,
respectively 5,  we observe that the values of the area relative to $\epsilon$ are 4, respectively 8 times the value
relative to $\frac{\epsilon}{2}$. This gives $\alpha=3$, respectively $\alpha=4$: the observed primitives appear numerically to be powers of the Dirac delta measure. Therefore we have numerically
put in evidence approximate solutions which appear to be derivatives of powers 2, 3 and 4 of the Dirac delta measure. They  are  obtained
as continuations of the Cauchy-Kovalevska solutions (since the same results can be obtained starting from an
analytic solution). In absence of definition such as the one in \cite {Egorov} as germs in $\epsilon$ of functions $f(x,\epsilon)=\delta(x,\epsilon)^n$, if $\delta(x,\epsilon)$ represents as usual approximations of the Dirac delta distribution, or, equivalently,   asymptotic objects such as in \cite{MaslovOmel},   these derivatives of powers of Dirac measures are likewise to be confused with numerical
blow up  since they reach very high top values and are not familiar mathematical objects, since they are not defined within distribution theory. Numerically, their primitives have the aspect  of Dirac delta distributions with  very high peaks and  an area between the peak and the $x$-axis that grows when the space step size  $\epsilon$ diminishes, as observed in the array. The fact that these new objects are continuations of classical analytic solutions after the analytic blow up is justified by the proof in section 6 below.\\
 
\textbf{4. Extension to Panov-Shelkovich $\delta^{(p)}$-shock waves.} A Panov-Shelkovich system for
$\delta^{"}$-shock waves can be stated, \cite{Panov} p. 83,
\begin{equation}\frac{\partial}{\partial t}u+\frac{\partial}{\partial x}(u^2)=0,  \end{equation}
\begin{equation}\frac{\partial}{\partial t}v +2\frac{\partial}{\partial x}(uv)=0,  \end{equation}
\begin{equation}\frac{\partial}{\partial t} w+2\frac{\partial}{\partial x}(uw)+2\frac{\partial}{\partial x}(v^n)=0, \ n=2, \end{equation}
\begin{equation} \frac{\partial}{\partial t}Z+2\frac{\partial}{\partial x}(uZ)+6\frac{\partial}{\partial x}(vw)=0 .\end{equation}

To obtain a weak asymptotic method it suffices to construct $u^\epsilon,v^\epsilon,w^\epsilon$ as in section 2,
then define $Z^\epsilon$ as the solution of the linear equation with second member

\begin{equation} \frac{d}{dt}Z^\epsilon(x,t)=\frac{2}{\epsilon}[(Z^\epsilon u^{\epsilon
+})(x-\epsilon,t)-(Z^\epsilon |u^{\epsilon }|)(x,t)+(Z^{\epsilon} u^{\epsilon
-})(x+\epsilon,t)]-6\frac{\partial}{\partial x}[(v^\epsilon
w^\epsilon)(.,t)*\phi_{\epsilon^\gamma}](x)\end{equation}
\\
for $\gamma>0 $ small enough, and with regularized initial condition $Z_0^\epsilon$. Of course as in section 2
the numerical coefficients play no role in the proof and the result extends without any change in proof to far
more general situations in one space dimension. In the proof below we set all coefficients equal one for
simplification. The coefficients in (22) will be used in the numerical test depicted in figure 2. We state same
assumptions on $u_0,v_0,w_0$ and their regularizations as above; $Z_0$ and its regularizations $Z_0^\epsilon$
are assumed to be $L^1_{loc}$ and periodic with same period. By induction the result holds clearly for
$\delta^{(p)}$ Panov-Shelkovich shock-waves with arbitrary $p\in\mathbb{N}$, \cite{Panov}. \\

\textbf{Theorem 2. Provided $\alpha,\beta, \gamma>0$ small enough and $\beta<\alpha$ the family
$(u^\epsilon,v^\epsilon,w^\epsilon,Z^{\epsilon})_\epsilon$ provides a weak asymptotic method for system
(18-21).}\\
\\
 \textit{proof.} We have to prove that $\forall \psi\in \mathcal{C}_c^\infty(\mathbb{R})$
\begin{equation} I:=\int\{\frac{\partial}{\partial t}(Z^\epsilon)\psi-(u^\epsilon Z^\epsilon)\psi'-(v^\epsilon w^\epsilon)\psi'\}
dx\rightarrow 0 \end{equation}
 when $\epsilon\rightarrow 0.$
  From (22), with  coefficients 2 and 6 replaced by 1,\\
 
$I=\int\{\frac{1}{\epsilon}[(Z^\epsilon u^{\epsilon +})(x-\epsilon,t)-(Z^\epsilon (u^{\epsilon +}+u^{\epsilon
-}))(x,t)+(Z^\epsilon u^{\epsilon -})(x+\epsilon,t)]\psi(x)+((v^\epsilon
w^\epsilon)*\phi_{\epsilon^\gamma})(x,t)\psi'(x)-(u^\epsilon Z^\epsilon)(x,t)\psi'(x) -(v^\epsilon
w^\epsilon)(x,t)\psi'(x)\}dx.$\\
\\
One can share $I$ into $I=I_1+I_2$ with \\
$$I_1=\int\{(Z^\epsilon u^{\epsilon +})(x,t)\frac{\psi(x+\epsilon)-\psi(x)}{\epsilon}-(Z^\epsilon u^{\epsilon
-})(x,t)\frac{\psi(x)-\psi(x-\epsilon)}{\epsilon}-(u^\epsilon Z^\epsilon)(x,t)\psi'(x)\}dx,$$

$$I_2=\int\{(v^\epsilon w^\epsilon)(x-y,t)\frac{1}{\epsilon^\gamma}\phi(\frac{y}{\epsilon^\gamma})-(v^\epsilon
w^\epsilon)(x,t)\frac{1}{\epsilon^\gamma}\phi(\frac{y}{\epsilon^\gamma})\}\psi'(x)dxdy.$$
\\
We first consider $I_2$. After two standard changes of variables\\
$$ I_2=\int\{(v^\epsilon w^\epsilon)(x,t)\phi(\mu)(\psi'(x+\epsilon^\gamma \mu)-\psi'(x))\}dxd\mu.$$
Therefore, from (12) and (16),\\
$$ |I_2|\leq const.\frac{1}{\epsilon^\alpha}.\epsilon^\gamma \int_{compact}|w^\epsilon(x,t)|dx \leq
const\frac{1}{\epsilon^\alpha}.\epsilon^\gamma \frac{1}{\epsilon^{(n+1)\alpha}}$$
i.e. 
 \begin{equation} I_2\leq const .\epsilon^{\gamma-(n+2)\alpha} \end{equation}
 which tends to 0 when $\epsilon\rightarrow 0$ provided $\gamma>(n+2)\alpha$. Now let us consider $I_1$.\\
 
$I_1=\int\{(Z^\epsilon u^{\epsilon+})(x,t)\psi'(x)+(Z^\epsilon u^{\epsilon +})(x,t)O_{1,[x]}(\epsilon)-(Z^\epsilon
u^{\epsilon-})(x,t)\psi'(x)+(Z^\epsilon u^{\epsilon -})(x,t)O_{2,[x]}(\epsilon)-(Z^\epsilon
u^\epsilon)(x,t)\psi'(x)\}dx.$\\
 \\
 After simplification, from (4, 7) 
 
 \begin{equation} |I_1|\leq \epsilon. const.\int|Z^\epsilon(x,t)|dx. \end{equation}
 
Now we need to evaluate $\int|Z^\epsilon(x,t)|dx$. To this end  we do as in lemma 1 and lemma 3: from (22) (dropping
the  coefficients for simplification)\\
 \\
$Z^\epsilon(x,t+dt)=Z^\epsilon(x,t)+\frac{dt}{\epsilon}[(Z^\epsilon u^{\epsilon +})(x-\epsilon,t)-(Z^\epsilon
|u^{\epsilon }|)(x,t)+(Z^\epsilon u^{\epsilon -})(x+\epsilon,t)] -dt\int((v^\epsilon
w^\epsilon)(x-y,t)\frac{1}{\epsilon^{2\gamma}}\phi'(\frac{y}{\epsilon^\gamma})dy+dt.o^\epsilon(x,t)(dt).$\\
 \\
The proof of this formula, where  $o^\epsilon(x,t)(dt)$ converges to 0 when $dt\rightarrow 0$ uniformly in $x\in \mathbb{R}$ and $t $ in a compact set, follows from the mean value theorem as exposed in the proof of lemma 1. Taking into account
simplifications in the integrals due to the periodicity, and for $dt>0$ small enough depending on $\epsilon$, so
that $1-\frac{dt}{\epsilon}\|u^\epsilon\|_\infty \geq 0$, one obtains as in lemmas 1, 3\\
  
$ \int_{-\pi}^{+\pi}|Z^\epsilon(x,t+dt)|dx \leq \int_{-\pi}^{+\pi}|Z^\epsilon(x,t)|dx +dt \|v^\epsilon\|_\infty
\int_{compact}|w^\epsilon(x-y,t)|dy \frac{const}{\epsilon^{2\gamma}}+dt. o^\epsilon (dt),$\\
 i.e. from (12) and (16)\\
 
$ \int_{-\pi}^{+\pi}|Z^\epsilon(x,t+dt)|dx \leq \int_{-\pi}^{+\pi}|Z^\epsilon(x,t)|dx
+dt\frac{const}{\epsilon^\alpha}\frac{1}{\epsilon^{(n+1)\alpha}} \frac{1}{\epsilon^{2\gamma}}+dt. o^\epsilon
(dt).$\\
 \vskip 4 cm

\textit{figure 2. Emergence of a $\delta^{"}$-wave from an analytic solution.}\\
\\

Finally as in proofs of lemmas 1, 3, sharing the interval $[t,t+\tau]$ into intervals of length $\frac{\tau}{n}$
with $n\rightarrow \infty$ one obtains\\
 
$ \int_{-\pi}^{+\pi}|Z^\epsilon(x,t+\tau)|dx \leq \int_{-\pi}^{+\pi}|Z^\epsilon(x,t)|dx +\tau
\frac{const}{\epsilon^{2\gamma+(n+2)\alpha}}.$\\
  \\    
 Replacing $t$ by $0$ and $\tau $ by $t$ one obtains  that\\
 
$ \int_{-\pi}^{+\pi}|Z^\epsilon(x,t)|dx \leq \int_{-\pi}^{+\pi}|Z_0^\epsilon(x)|dx + t
\frac{const}{\epsilon^{2\gamma+(n+2)\alpha}}.$\\
  \\
  From (25)
 \begin{equation} |I_1|\leq const.\epsilon+const. \epsilon^{1-2\gamma-(n+2)\alpha} \end{equation}
when $t$ ranges in a bounded interval. From (24) and (26) $I\rightarrow 0$ when $\epsilon\rightarrow 0$
provided $\gamma>(n+2)\alpha, 2\gamma+(n+2)\alpha <1$ and, further,  $\beta<\alpha$ to apply lemma 2.$\Box$
 \begin{figure}[h]
\centering \includepdf[width=\textwidth]{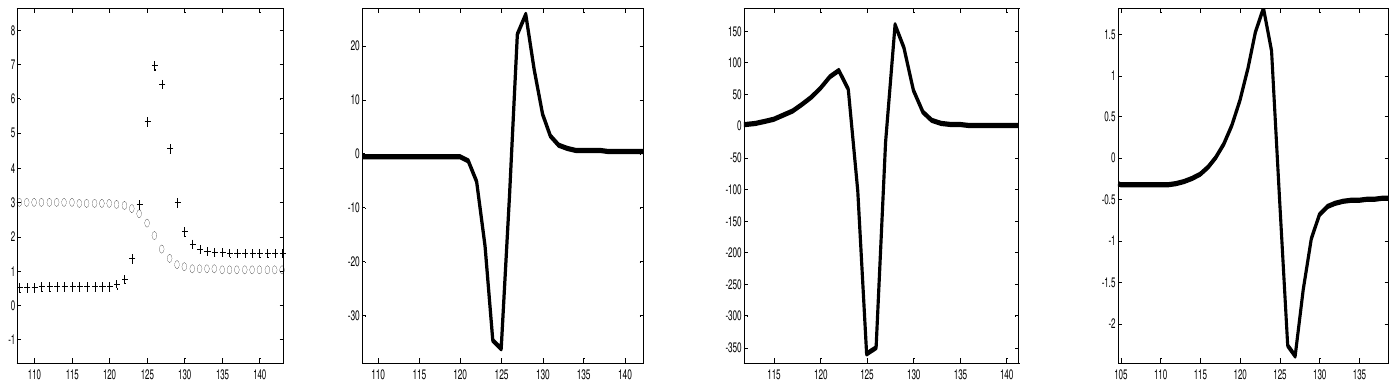}
\end{figure}

\textit{Numerical confirmation.} In figure 2 we start from an analytic solution $(u,v,w,Z)$ of system (18-21) and we
observe the creation of the well known shock wave in $u$ and of a $\delta$-wave in $v$ (left panel), of a
$\delta'$-wave in $w$ (left-middle panel) and of a $\delta^{"}$-wave in Z (right-middle panel). If $\delta$ and
$\delta'$-waves are very clearly identified, the observation of $\delta^{"}$-waves (right-middle panel) is not
always very clear: one observes two small positive peaks separated by a deep negative peak because
$\delta^{"}$-waves are derivatives of $\delta'$-waves whose graph is nearly vertical in its middle (left-middle
panel): the very large value of the derivative there can relatively hide the two positive peaks of a
$\delta^{"}$-wave. So an idea to observe them better is to compute a primitive of $Z$ (right panel) to check
that it has a $\delta'$-shape. It is also convenient to observe a primitive of a primitive of $Z$: one obtains
an affine function with  a delta peak located on the shock wave.\\

\textbf {5.  Multi-dimensional systems and other extensions.} One can extend system (1, 2, 3) to 2-D as 
\begin{equation} \frac{\partial}{\partial t}u +\frac{\partial}{\partial x}(f_1(u))+\frac{\partial}{\partial y}(f_2(u))=0, \ \frac{\partial}{\partial t} v +\frac{\partial}{\partial x}(g_1(v))+\frac{\partial}{\partial y}(g_2(v))=0, \end{equation}

\begin{equation} \frac{\partial}{\partial t} \rho+\frac{\partial}{\partial x}(\rho u) +\frac{\partial}{\partial y}(\rho v)=0,\end{equation} 
\begin{equation}\frac{\partial}{\partial t} w +\frac{\partial}{\partial x}(uw)+ \frac{\partial}{\partial y}(vw) +\frac{\partial}{\partial x}(P(\rho))+\frac{\partial}{\partial y}(Q(\rho)) =0, \end{equation}
\\
where $(u,v)$ plays the role of $u$ in system (1, 2, 3), $f_1,f_2,g_1$ and $g_2$ are smooth functions ; $\rho$ plays
the role of $v$ ; $P$ and $Q$ are polynomials. The equations (27) give a 2-D discontinuous "`velocity"' $(u,v)$, (28) is a
continuity equation producing delta-waves and (29) is an equation extending (3) to 2-D. We seek a
solution on the 2-dimensional torus  $\mathbb{T}^2=\mathbb{R}^2/(2\pi\mathbb{Z})^2$. \\

We use solutions of the scalar conservation laws (27) having the properties (4, 5) in 2-D: uniform boundedness
and
\begin{equation} \forall \beta>0 \ \forall \delta>0 \ \exists \ const / |\frac{\partial}{\partial
x}u^\epsilon(x,y,t)|\leq \frac{const}{\epsilon^\beta}, \ \ |\frac{\partial}{\partial y}u^\epsilon(x,y,t)|\leq
\frac{const}{\epsilon^\beta} \ \ \ \forall (x,y)\in \mathbb{T}^2 \ \forall t\in [0,\delta],\end{equation}
with  same property for $v^\epsilon$. For equation (28) we state the ODE as follows (\cite{Colombeaugravitation} section 7), which consists in stating (8) in both axis directions :\\ 

$\frac{d}{dt}X^\epsilon(x,y,t)=\frac{1}{\epsilon}[X^\epsilon(x-\epsilon,y,t)u^{\epsilon+}(x-\epsilon,y,t)-X^\epsilon(x,y,t)|u^\epsilon(x,y,t)|$\\
$+X^\epsilon(x+\epsilon,y,t)u^{\epsilon-}(x+\epsilon,y,t)+X^\epsilon(x,y-\epsilon,t)v^{\epsilon+}(x,y-\epsilon,t)-X^\epsilon(x,y,t)|v^{\epsilon}(x,y,t)|+X^\epsilon(x,y+\epsilon,t)v^{\epsilon-}(x,y+\epsilon,t)]$,\begin{equation}
\end{equation}

\begin{equation} X^\epsilon(x,y,0)=\rho_0^\epsilon(x,y)\end{equation}
and
\begin{equation} \rho^\epsilon(x,y,t)=(X^\epsilon(.,.,t)*\phi_{\epsilon^\alpha})(x,y)\end{equation}
where $\phi$ is a $\mathcal{C}^\infty $ function on $\mathbb{R}^2 $ with compact support such that
$\int\phi(x,y)dxdy=1,
\phi_{\epsilon^\alpha}(x,y)=\frac{1}{\epsilon^{2\alpha}}\phi(\frac{x}{\epsilon^\alpha},\frac{y}{\epsilon^\alpha})$.\\

We  sketch the proofs which are direct extensions of the proofs in the 1-D case.\\

\textbf{Lemma 4.} $\int_{-\pi}^{+\pi}\int_{-\pi}^{+\pi}|X^\epsilon(x,y,t)|dxdy\leq
\int_{-\pi}^{+\pi}\int_{-\pi}^{+\pi}|\rho_0^\epsilon(x,y)|dxdy\leq constant $\\ $independent \ on \ \epsilon.$\\
\\
proof. Once one has written the ODE in the form\\

$X^\epsilon(x,y,t+dt)=\frac{dt}{\epsilon}X^\epsilon u^{\epsilon+}(x-\epsilon,y,t)+\frac{dt}{\epsilon}X^\epsilon
v^{\epsilon+}(x,y-\epsilon,t)+
(1-\frac{dt}{\epsilon}|u^\epsilon(x,y,t)|-\frac{dt}{\epsilon}|v^\epsilon(x,y,t)|)X^\epsilon(x,y,t)+\frac{dt}{\epsilon}X^\epsilon
u^{\epsilon-}(x+\epsilon,y,t)+\frac{dt}{\epsilon}X^\epsilon
v^{\epsilon-}(x,y+\epsilon,t)+dt.o^\epsilon(x,y,t)(dt)$,\begin{equation}\end{equation}
 the proof, based on periodicity, is the same as the proof of lemma 1.$\Box$\\
 \\
Corollary. $|\rho^\epsilon(x,y,t)|\leq \frac{M}{\epsilon^{2\alpha}}, \ |\frac{\partial}{\partial
x}\rho^\epsilon(x,y,t)|, \ |\frac{\partial}{\partial y}\rho^\epsilon(x,y,t)|\leq
\frac{M}{\epsilon^{3\alpha}}.$\\

The proof follows from the convolution (33). Using the proof of lemma 2 in the 2 variables $x,y$ one obtains
that the family $(u^\epsilon,v^\epsilon,\rho^\epsilon)$ is a weak asymptotic method for equation (28).\\

For equation (29) we state\\
$ \frac{d}{dt}w^\epsilon(x,y,t)=\frac{1}{\epsilon}[(w^\epsilon u^{\epsilon +})(x-\epsilon,y,t)- (w^\epsilon
|u^{\epsilon}|)(x,y,t)+(w^\epsilon u^{\epsilon -})(x+\epsilon,y,t)+(w^\epsilon v^{\epsilon +})(x,y-\epsilon,t)-
(w^\epsilon |v^{\epsilon}|)(x,y,t)+(w^\epsilon v^{\epsilon -})(x,y+\epsilon,t)]- P'(\rho^\epsilon(x,y,t))
\frac{\partial}{\partial x}\rho^\epsilon(x,y,t)-$
\begin{equation}Q'(\rho^\epsilon(x,y,t)) \frac{\partial}{\partial y}\rho^\epsilon(x,y,t).\end{equation}

\textbf{Lemma 5.} $ \exists N\in \mathbb{N}, \exists C>0 \ / \ \forall \delta>0, \forall t\in[0,\delta], \forall \epsilon>0 \ \ \int_{-\pi}^{+\pi}|w^\epsilon(x,y,t)|dx \leq \frac{C}{\epsilon^{N\alpha}}$.\\

It suffices to develop (35) according to Taylor's formula, as done in  (34), with the two additional terms involving $P$ and $Q$ and to follow the
proof of lemma 3.$\Box$\\

Then a proof similar to that of theorem 1 shows that  (35) provides a weak asymptotic method to equation (29).\\

The weak asymptotic method under consideration applies in arbitrary dimension on $\mathbb{T}^d$ and $\mathbb{R}^d$. For simplicity and to facilitate the reading we consider a $2\times 2$ system on the 1-D torus since multidimensional extensions are rather easy \cite{Colombeaugravitation}.
 Now we show that the weak asymptotic method  extends easily to  systems of the form 

\begin{equation} \frac{\partial}{\partial t}u+\frac{\partial}{\partial x}(uf(u,v))=0\end{equation}
\begin{equation}\frac{\partial}{\partial t} v+\frac{\partial}{\partial x}(vg(u,v))=0,\end{equation}
on the 1-D torus, with $f,g$ analytic bounded on $\mathbb{R}^2$, since this assumption  will permit a very simple proof. This weak asymptotic method has been adapted to pressureless fluids in \cite{Colombeaugravitation}. For system (36, 37) we consider the ODE

\begin{equation}\frac{du^\epsilon}{dt}(x,t)=\frac{1}{\epsilon}[u^\epsilon(x-\epsilon,t)f^{\epsilon,+}(x-\epsilon,t)-u^\epsilon(x,t)|f^{\epsilon}(x,t)|+u^\epsilon(x+\epsilon,t)f^{\epsilon,-}(x+\epsilon,t)],\end{equation}
where $f^\epsilon=f(u^\epsilon,v^\epsilon)$ and similar ODE for (37). We sketch the proof: a detailed proof for the (more difficult) system of pressureless fluids is given in \cite{Colombeaugravitation}. The initial data $u_0^\epsilon, v_0^\epsilon\in \mathcal{C}_b(\mathbb{T})$  with $\exists  C>0 \  / \|u_0^\epsilon\|_{L^1(\mathbb{T})},\|v_0^\epsilon\|_{L^1(\mathbb{T})} < C \forall \epsilon>0$.\\

For fixed $\epsilon>0$ the system of two ODEs (38) admits a local solution $u^\epsilon,v^\epsilon: [0,\delta[\longmapsto \mathcal{C}_b(\mathbb{R})\times \mathcal{C}_b(\mathbb{R})$.  First we prove that
\begin{equation}\|u^\epsilon(.,t)\|_\infty\leq \|u_0^\epsilon\|_\infty exp(\frac{const}{\epsilon}t).\end{equation}
Indeed, from (38) $|u^\epsilon(x,t)|\leq |u_0^\epsilon(x)|+\frac{const}{\epsilon}\int_0^t\|u^\epsilon(.,s)\|_\infty ds$ since $f$ is assumed bounded; then one applies the Gronwall formula. It follows easily from (39) that (38) has a global solution on $[0,+\infty[$, as in \cite{Colombeaugravitation} section 4. Now we prove that 

 \begin{equation}\int_0^{2\pi}|u^\epsilon(x,t)|dx\leq \int_0^{2\pi}|u_0^\epsilon(x)|dx.\end{equation}
\\
Indeed, from (38) the proof of (40) is identical to the proof of lemma 1. It follows  that the family $(u^\epsilon,v^\epsilon)$ is a weak asymptotic method for system (36, 37): one has to prove that 
$$\forall \psi \in \mathcal{C}_c^\infty(\mathbb{R}) \ \int (\frac{\partial}{\partial t}(u^\epsilon) \psi-u^\epsilon f^\epsilon \psi') dx\rightarrow 0 $$ when $\epsilon \rightarrow 0.$
From (38, 40) the proof is similar to that of lemma 2.\\

In \cite{ColombeauSiam, ColombeauNMPDE} a numerical scheme for the 3-D system of pressureless fluids and the 3-D Euler-Poisson system has been investigated mathematically and numerically. In the case of 3-D pressureless fluids without selfgravitation and 1-D Euler-Poisson equations it is proved there that this scheme provides a weak asymptotic method with weak derivative both in space and time: th. 3 p. 1909 in \cite{ColombeauSiam} with annoucement p. 1911 of the result in 3-D proved in  \cite{ColombeauNMPDE} p. 96-100, th. 2 p. 86 in \cite{ColombeauNMPDE} for the 1-D Euler-Poisson equations. Since the author was unaware of the concept of weak asymptotic method \cite{Danilov1} the weak asymptotic method there is called "convergence" in \cite{ColombeauSiam} and "consistence" in \cite{ColombeauNMPDE}.\\

For the 2-D and 3-D Euler-Poisson system the scheme in \cite{ColombeauNMPDE} fails to provide a weak asymptotic method since one is forced (th. 1 p. 85) to assume the boundedness of the velocity vector and of the gradient of the gravitation potential. This problem has been solved in \cite{Colombeaugravitation} by an adaptation of the method presented in this paper. This adaptation relies on a priori estimates which have been avoided in the example (36, 37) by the assumption that $f$ and $g$ are bounded on $\mathbb{R}^2$.\\

\textbf{6. The weak asymptotic method in the analytic case.} In section 2 we have proved that the solutions of the system of ODEs is a weak asymptotic method for the system (1-3), which has received a confirmation in section 3 from a numerical solution of the ODEs of section 2. Now we consider the classical case and we prove that the weak asymptotic method under consideration gives the classical analytic solution when it exists. Indeed the numerical scheme of section 3 gives a smooth solution; now we prove this smooth solution  is  the analytic solution as it should be.  The space of all holomorphic maps from $\Omega\subset\mathbb{C}$ open into a Banach space $E$ is denoted by $\mathcal{H}(\Omega,E)$.  We assume holomorphy in the initial conditions and coefficient:\\

$v_0^\epsilon, w_0^\epsilon\in  \mathcal{H}(\mathbb{T}\times\{y\in ]-r,+r[\})$ and  $u^\epsilon(.,.)\in\mathcal{H}(\mathbb{T}\times \{y\in ]-r,+r[\}\times \{t/ |t|<a\})$ for some $r,a>0$ uniformly in $\epsilon$.\\

 We assume the coefficients have a fixed sign  away from 0:$$\exists \eta>0 \ / \  |u^\epsilon(x,t)|\geq \eta \  \forall (x,t)\in \mathbb{T}\times\{t/ |t|<a\},$$ so as to get rid of the lack of
analyticity due to $u^\pm$ and $|u|$ in (8, 14): in the following proof we assume they are positive, which is the case in figures 1 and 2. Then
equation (8) becomes
\begin{equation} \frac{\partial X^\epsilon}{\partial
t}(x,t)=-\frac{1}{\epsilon}[X^\epsilon(x,t)u^\epsilon(x,t)-X^\epsilon(x-\epsilon,t)u^\epsilon(x-\epsilon,t)].\end{equation}
It has previously been proved that a global solution exists in the space $\mathcal{C}_b(\mathbb{R})$. But we do
not know the existence and the nature of a limit when $\epsilon\rightarrow 0$. The purpose of this section is
to prove that \textit{the approximate solutions converge to the classical analytic solution}. This is
done by applying an abstract version of the Cauchy-Kovalevska theorem. Let us recall the version of this theorem presented in \cite{Treves} theorem 17.2 p. 148.\\

 \textbf{Definition.} A scale of Banach spaces is a family of Banach spaces $(E_s)_s, 0<s<s_0$, such that $\forall s,s'\in ]0,s_0], \ s>s' \Rightarrow E_s\subset E_{s'}$ with inclusion $i_{s,s'}: E_s\longmapsto  E_{s'}$ such that $\|i_{s,s'}\|_{L(E_s,E_{s'})}\leq 1$ where $L(E_s,E_{s'})$ is the Banach space of all linear continuous maps from $E_s$ into $E_{s'}$.\\

\textbf{Theorem:  Abstract Cauchy-Kovalevska.}  We consider the Cauchy problem
 \begin{equation} u'(t)=A(t)u(t)+ f(t), \ u(0)=u_0\in E_{s_0}, \end{equation}
where, for some $a>0$,  $\forall t\in \{|t|<a\}\subset \mathbb{C}, \ A(t)\in  L(E_s,E_{s'})$ as soon as $s<s'$, and the map $A\in    \mathcal{H}(\{|t|<a\}, L(E_s,E_{s'})).$ The map $f\in \mathcal{H}(\{|t|<a\},E_{s_0})$. The main assumption is:\\
 \begin{equation}\exists M>0 \ / \ \forall t\in \{|t|<a\}, \  \forall s>s' \ \ \|A(t)\|_{ L(E_s,E_{s'})}\leq \frac{M}{s-s'}.\end{equation}

Then $\exists \  C>0$ depending only on $M$ (not on $f$ and $u_0$)  such that $\forall s<s_0$ if $\delta:= min(a,\frac{C}{s_0-s}) \ \ \exists !$    solution $u\in \mathcal{H}(\{|t|<\delta\},E_{s})$. Further the bounds of $u$ depend only on domains and bounds of the data.\\

To apply this theorem we
consider the scale of Banach spaces $(E_s)_{0<s<s_0}$, for some $s_0>0$, defined by

 \begin{equation} E_s:=\{f\in  \mathcal{H}(\mathbb{T}\times]-s,s[,\mathbb{C}) \  
  continuous \ and \  bounded \ on \ the \ closure \ of \ this \ strip\}, \end{equation}
\\
equipped with the sup norm on the strip. The real number $s_0>0$ is chosen so that $\forall \epsilon, t$ the functions $(u^\epsilon(.,t),
v_0^\epsilon, w_0^\epsilon)_\epsilon$ are elements of the space $E_{s_0}$, uniformly bounded in sup norm. Then, if $ s>s'$,  for all $t$, the map $\mathcal{A}^\epsilon(t) \in  L(E_s,E_{s'})$ is defined by:
\begin{equation} [z\mapsto X(z)]\in E_s\longmapsto [z\mapsto
-\frac{1}{\epsilon}(X(z)u^\epsilon(z,t)-X(z-\epsilon)u^\epsilon(z-\epsilon,t))]\in E_{s'}.
\end{equation} 
 We have $ \|\mathcal{A}^\epsilon(t)\|_{L(E_s,E_{s'})}\leq \frac{const}{s-s'}$, where $const$ is
independent on $s,s',t$ and $\epsilon$. This follows at once from the mean value theorem applied to the second member of (45) and  from Cauchy's inequality for the derivative of a
holomorphic function.\\

 The abstract version
of the Cauchy-Kovalevska theorem  gives : \\ $\exists C>0, \ / \ \forall \epsilon, \ \forall   s<s_0  \  \exists ! X^\epsilon\in\mathcal{H}(\{|t|<\frac{C}{s_0-s}\},E_{s})$
 solution of (41) with initial condition
$v_0^\epsilon$, where   $C$ is independent on $s$ and $\epsilon$. Setting $\delta=\frac{C}{s_0-s}$, the functions $X^\epsilon: (z,t)\longmapsto X^\epsilon(z,t)$ so obtained
 are holomorphic  on  $\mathbb{T}\times]-s,s[\times\{|t|<\delta\} $with values in
$\mathbb{C}$, bounded on compact sets uniformly in $\epsilon$, from the uniform bounds and domains  of $v_0^\epsilon$ and $u^\epsilon(.,t)$ in $\epsilon$ and $t$. From the theory of normal families of holomorphic functions, from
any sequence of the $X^\epsilon$s one can extract a subsequence that converges uniformly on the compact sets of  $\mathbb{T}\times]-s,s[\times\{|t|<\delta\}. $
 The limit when $\epsilon\rightarrow 0$ is analytic and unique extension to the complex domain of the classical solution of
the linear Cauchy-Kovalevska problem
\begin{equation} \frac{\partial X}{\partial t}(x,t)=-\frac{\partial (X u)}{\partial x}(x,t), \
X(x,0)=v_0(x).\end{equation}
Therefore all subsequences of $(X^\epsilon) $ converge to the same limit. Therefore the whole sequence
$(X^\epsilon)$ converges uniformly on compact sets of   $\mathbb{T}\times\{0<t<\delta\} $ to the solution of (46) when $\epsilon\rightarrow 0$, i.e. of (2), with $b=1$ for simplicity. The proof
applies as well to (14): for the second member of (14) one uses the uniform bounds and uniform domain on the $v^\epsilon$s obtained from the abstract Cauchy-Kovalevska theorem applied to (8). The  proof ensures the same existence time for the weak asymptotic method as when it is applied to prove existence of the classical analytic solution. It also applies  in several space dimension. \\

The proof in this section can be reproduced practically without any change for systems  such as (36, 37) and those  in \cite{Colombeaugravitation},  using the nonlinear theorem due to L. Nirenberg and
T. Nishida  \cite{Nirenberg, Nishida} instead of the linear abstract Cauchy-Kovalevska theorem  \cite{Treves}.\\

\textit{In conclusion we have proved that with analytic initial data (and constant sign of $u$) the weak asymptotic method starts by giving the
classical analytic solution, which proves that the numerically observed very irregular shock-waves are
continuations of the classical analytic solutions after the analytic blow up occurs, as observed numerically.} \\


\textbf{7. Conclusion.} This new technique of construction of weak asymptotic methods permits to study PDEs having $\delta,\delta',\dots$ -shock wave solutions, and other systems, such as the 3-D Euler-Poisson system \cite{Colombeaugravitation} for applications in fluid dynamics and cosmology \cite{ColombeauNMPDE, Colombeauideal, Colombeaugravitation}, by transfering the problem to a family of ODEs in Banach spaces. This permits to demonstrate in a mathematically rigorous way the existence of approximate solutions for the Cauchy problem, extending previous results on the Riemann problem.\\

 Furthermore this  method permits to construct very simple numerical schemes that bring a confirmation of the theoretical results and permit a visualization of the approximate solutions. In particular it has permitted to put in evidence Panov-Shelkovich $\delta^{(n)}$-shocks and  much more  irregular shock waves (from a more general family of systems of conservation laws) as continuations of
the classical analytic Cauchy-Kovalevska solutions when their existence time is over. \\

Of course the problem of uniqueness of a \lq\lq{}limit\rq\rq{} of these approximate solutions remains unsolved. It has been checked numerically on standard systems of fluid dynamics that the adaptation of the method presented here to equations of fluid dynamics \cite{Colombeaugravitation} and to the very closely related numerical method in \cite{ColombeauSiam, ColombeauNMPDE, Colombeauideal} have always given the known solutions even on very demanding tests, such as tests of Woodward-Colella, Toro, Lax, see \cite{Colombeauideal}.\\


\begin{thebibliography}{<50 >}

\bibitem{AlbeDani} S. Albeverio, V.G. Danilov. Global in time solutions to Kolmogorov-Feller pseudodifferential
equations with small parameter. Arxiv: 1101.5836V1, 31 jan. 2011.

\bibitem{Albeverio} S. Albeverio, O.S. Rozanova, V.M. Shelkovich. Transport and concentration processes in the
multidimensional zero-pressure gas dynamics model with the energy conservation law. ArXiv: 1101.581v1, 30 Jan
2011.

\bibitem{AlbeShelk} S. Albeverio, V.M. Shelkovich. On delta shock front problem. In "`Analytical approach to
multibalance laws"', chap. 2. Editor O. Rozanova, Nova Science Pub. Inc. 2005, pp. 45-88.






\bibitem{ColombeauSiam} M. Colombeau. A method of projection of delta waves in a Godunov scheme and application
to pressureless fluid dynamics. SIAM J. Numer. Anal. 48, 5, 2010, pp. 1900-1919.

\bibitem{ColombeauNMPDE} M. Colombeau. A consistent numerical scheme for self-gravitating fluid dynamics. Num.
Methods for PDEs. 29, 1, 2013, pp. 79-101.


\bibitem{Colombeauideal} M. Colombeau. A simple numerical scheme for the 3-D system of ideal gases and a study
of approximation in the sense of distributions. J. Comput. Appli. Math. 248, 2013, pp. 15-30.
\bibitem{Colombeaugravitation} M. Colombeau. Weak asymptotic methods for 3-D self-gravitating pressureless
fluids. Application to the creation and evolution of solar systems from the fully nonlinear Euler-Poisson
equations. Preprint.





\bibitem{Mitrovic} V.G. Danilov, D. Mitrovic. Delta shock wave formation in the case of triangular hyperbolic
system of conservation laws. J. Differential Equations 245, 2008, pp. 3704-3734.

\bibitem{Mitrovic2} V.G. Danilov, D. Mitrovic. Shock wave formation process for a multidimensional scalar
conservation law. Quarterly of Applied Mathematics vol 69, 2011, pp. 613-634.



\bibitem{DanilovO1}  V. G. Danilov, G.A. Omel'yanov. Weak asymptotic method and the interaction of infinitely narrow $\delta$-solitons. Nonlinear Analysis 54, 2003, pp. 773-799.


\bibitem{DanilovO2}  V. G. Danilov, G.A. Omel'yanov.  Weak asymptotic method for the study of infinitely narrow $\delta$-solitons. Electronic J. Diff. Equations 90, 2003, 27 p.


\bibitem{Danilov1} V.G. Danilov, G.A. Omel'yanov, and V.M. Shelkovich. Weak Asymptotic Method and Interaction
of Nonlinear Waves, AMS Translations vol 208, 2003, pp. 33-164.




\bibitem{Shelkovich2} V.G. Danilov, V.M. Shelkovich. Dynamics of propagation and interaction of $ \delta$ shock
waves in conservation law systems. J. Differential Equations 211, 2005, pp. 333-381.

\bibitem{Shelkovich3} V.G. Danilov, V.M. Shelkovich. Delta-shock wave type solution of hyperbolic systems of
conservation laws. Quart. Appl. Math. 63, 2005, pp. 401-427.




\bibitem{Egorov}  Y. Egorov.  A Theory of Generalized Functions. Russian Math. Surveys 45, 5, 1990, pp. 1-49.

\bibitem{Maslov} V.P. Maslov. Asymptotic Methods and Perturbation Theory. Nauka, Moscow, 1988.

\bibitem{MaslovOmel} V.P. Maslov, G.A. Omel'yanov. Asymptotic soliton form solutions of equations with small dispersion. Russian Math. Surveys 36, 3, 1981, pp. 73-119.

\bibitem{Mitrovic3} D. Mitrovic, V. Bojkovic, V.G. Danilov. Linearization of the Riemann problem for a triangular system of conservation laws and delta shock wave formation process. Mathematical Methods in the Applied Sciences.  2010, 33, pp. 904-921.
\bibitem{Nirenberg} L. Nirenberg.  An abstract form of the nonlinear Cauchy-Kovalevski theorem. J. Diff. Geometry 6, 1972, pp. 561-576.
\bibitem{Nishida} T. Nishida. A note on a theorem of Nirenberg. J. Diff. Geometry. 12, 1977, pp. 629-633.
\bibitem{Omel'yanov} G.A. Omel'yanov, I. Segundo-Caballero. Asymptotic and numerical description of the
kink/antikink interaction. Electronic J. of Differential Equations, 2010, 150, pp. 1-19.

\bibitem{Panov} E.Yu. Panov, V.M. Shelkovich. $\delta$'-shock waves as a new type of solutions to systems of
conservation laws. J. Differential Equations 228, 2006, pp. 49-86.






\bibitem{ShelkovichRMS} V.M. Shelkovich. $\delta-$ and $\delta'-$shock wave types of singular solutions of
systems of conservation laws and transport and concentration processes. Russian Math. surveys 63,3, 2008, pp.
405-601.
\bibitem{Shelkovichmat} V.M. Shelkovich. The Riemann problem admitting $\delta-,\delta$'-shocks and vacuum
states; the vanishing viscosity approach. J. Differential Equations 231, 2006, pp. 459-500.


\bibitem{Shelkovich1} V.M. Shelkovich. Transport of mass, momentum and energy in zero-pressure gas dynamics. In
Proceedings of Symposia in Applied Mathematics 2009; vol.67. Hyperbolic Problems: Theory, Numerics and
Applications. Edited by E. Tadmor, Jian-Guo Liu, A.E. Tzavaras. AMS, 2009, pp. 929-938.
\bibitem{Treves} F. Treves. Basic Linear Partial Differential Equations. Academic Press, 1975.



























  
\end{thebibliography}
 \end{document}